\documentclass{article}
\usepackage{amsmath,amsfonts,amssymb}
\usepackage[round]{natbib}
\usepackage{booktabs}

\renewcommand{\ge}{\geqslant}
\renewcommand{\le}{\leqslant}

\newcommand{\dnorm}{\mathcal{N}}

\newcommand{\za}{z_\alpha}
\newcommand{\zoma}{z_{1-\alpha}}

\title{Coverage errors for Student's $t$ confidence intervals comparable
to those in Hall (1988)}
\date{January 2025}
\author{Art B. Owen\\Stanford University}
\begin{document}
\maketitle
\begin{abstract}
  Table 1 of \cite{hall:1988} contains asymptotic coverage error formulas
  for some nonparametric approximate 95\% confidence intervals for the mean
  based on $n$ IID samples. The table includes an entry for
  an interval based on the central limit theorem using Gaussian
  quantiles and the Gaussian maximum likelihood variance estimate.
  It is missing an entry  
  for the very widely used Student's $t$ confidence intervals. This
  note develops such a formula. 
The impetus to  revisit this issue arose from the surprisingly robust performance
  of confidence intervals based on 
  Student's t statistic in randomized quasi-Monte Carlo sampling.
Hall's table had  $0.14\kappa -2.12\gamma^2-3.35$ 
for normal theory intervals; the corresponding
entry for Student's $t$ is $0.14\kappa -2.12\gamma^2$.

  An earlier version of this note reported that it corrected
  some coverage error formulas in \cite{hall:1988}.
  Two-sided errors take the form
 $2\Phi^{-1}(0.975)(A\kappa + B\gamma^2+C)\varphi(1.96)/n +O(1/n^{3/2})$
where the error may well be $O(n^{-2})$.
Hall's table showed
 $\Phi^{-1}(0.975)(A\kappa + B\gamma^2+C)$.
The version intended as a correction 
had $2(A\kappa + B\gamma^2+C)$, wider by about $2/1.96\doteq1.02$.
So, Hall's table really is proportional to the two-sided coverage errors.
\end{abstract}

\section{Introduction}

Table 1 of \cite{hall:1988} gives expressions for the coverage error in several
nonparametric confidence interval methods for the mean of a real-valued
random variable $X$ based on an IID sample of $n$ realizations of that
random variable.  It is missing an entry for the usual confidence interval
based on Student's $t$ statistic.  It does have an entry for the case where
the $t$ threshold is replaced by one from the standard
Gaussian distribution and the usual unbiased sample variance is replaced
by the maximum likelihood estimate (MLE). This note uses an expression from page 949 of that paper to develop a formula for the Student $t$ intervals.

The motivation to revisit this problem comes from recent work in randomized
quasi-Monte Carlo (RQMC) sampling in \cite{ci4rqmc}. RQMC provides an
unbiased estimate of a multidimensional integral. The sample values that go
into that estimate have a complex dependence structure that improves the
convergence rate beyond that of plain Monte Carlo. Then a small number
of independent replicates are made in order to estimate error.  The surprise
in~\cite{ci4rqmc} was that a plain Student's $t$ confidence interval proved
to be remarkably robust in terms of coverage accuracy at 95\% nominal coverage.
It performed better than two bootstrap confidence intervals.  A theoretical
comparison of bootstrap confidence intervals from \cite{hall:1988} revolves
around the skewness and kurtosis of the random variable being replicated.
This paper fills in the missing formula for Student's $t$ distribution.

For $n\ge2$, let $X_1,\dots,X_n$ be IID random variables with the same distribution as $X$.
Then for $0<\alpha<1/2$,  a standard Student's $t$ based confidence interval for $\mu$ 
with nominal coverage $1-2\alpha$ takes the form
\begin{align}\label{eq:t}
\bar X \pm t^{1-\alpha}_{(n-1)} s/\sqrt{n}
\end{align}
where
$$
\bar X = \frac1n\sum_{i=1}^nX_i\quad\text{and}\quad s^2=\frac1{n-1}\sum_{i=1}^n(X_i-\bar X)^2
$$
and $t^{1-\alpha}_{(n-1)}$ is the $1-\alpha$ quantile of Student's $t$ distribution on $n-1$
degrees of freedom.  This interval has asymptotic coverage probability $1-2\alpha$.
It is more common to write about coverage $1-\alpha$ using $t_{(n-1)}^{1-\alpha/2}$ for $\alpha\in(0,1)$
but the usage above conforms to Table 1 of \cite{hall:1988}.

An alternative interval is
\begin{align}\label{eq:norm}
\bar X \pm \Phi^{-1}(1-\alpha)\hat\sigma/\sqrt{n}
\end{align}
where $\Phi$ is the cumulative distribution function of the $\dnorm(0,1)$ distribution,
and
$$
\hat\sigma^2 = \frac1n\sum_{i=1}^n(X_i-\bar X)^2
$$
is the MLE of $\sigma^2$ in the model $X_i\sim\dnorm(\mu,\sigma^2)$.
This is the interval that appears in Table 1 of \cite{hall:1988} with the label `Norm'.

Table 1 of \cite{hall:1988} includes asymptotic coverage error formulas for `Norm'
and some bootstrap confidence intervals.
For the necessary assumptions on the distribution of $X$, see \cite{hall:1988}.
These table entries are for $\alpha=0.025$ corresponding to the common default
of 95\% confidence.  The coverage error of~\eqref{eq:norm} for $\alpha=0.025$ is
\begin{align}\label{eq:coverror}
\frac2n\Bigl(0.14\kappa -2.12\gamma^2-3.35\Bigr)\varphi(1.96) + O(n^{-2}) ,
\end{align}
where $\varphi$ is the $\dnorm(0,1)$ probability density function,
and the constants have been rounded.
The row labeled `Norm' in Table 1 of \cite{hall:1988}
contains $0.14\kappa -2.12\gamma^2-3.35$,
omitting the factor $2\varphi(1.96)/n$ common to all of the
coverage error expressions.

Table 1 is missing a row for the usual Student's $t$ interval in equation~\eqref{eq:t}.
That interval differs from the one for~\eqref{eq:norm} in using $s$
instead of $\hat\sigma$ and in using $t^{1-\alpha}_{(n-1)}$ instead of $\Phi^{-1}(1-\alpha)$.

Section~\ref{sec:deriv} rederives Hall's formula
for coverage errors in 95\% confidence intervals based on the Gaussian interval.
Section~\ref{sec:adjusted} then constructs formulas for the Student $t$
interval.  The key step is to use an asymptotic formula from \cite{abra:steg:1972}
relating quantiles of Student's $t$ distribution to those of the $\dnorm(0,1)$
distribution.  Both tasks were made simpler by noticing that Hall's derivation
for the `Norm' interval does not use all of the machinery that the bootstrap
intervals in that paper require.

\section{Derivation}\label{sec:deriv}

Here we look at the derivation for the row called `Norm' before
constructing the analogous formula for the standard $t$ intervals.
Let $z_\alpha=\Phi^{-1}(\alpha)$ for $0<\alpha<1$
and then let $C(\alpha) = \Pr(\mu\le\hat\theta + n^{-1/2}\hat\sigma z_{\alpha})$
be the one sided coverage probability for the normal interval
$(-\infty,\hat\theta + n^{-1/2}\hat\sigma z_{\alpha}]$.
From page 949 of \cite{hall:1988}, 
\begin{multline}\label{eq:hall949}
C(\alpha)=\alpha - n^{-1/2}\frac16\gamma(2z_\alpha^2+1)\varphi(z_\alpha) 
+n^{-1}z_\alpha\Bigl(\frac1{12}\kappa(z_\alpha^2-3) \\
-\frac1{18}\gamma^2\bigl( z_\alpha^4+2z_\alpha^2-3\bigr) 
-\frac14\bigl( z_\alpha^2+3\bigr)\Bigr)\varphi(z_\alpha)+O(n^{-3/2}). 
\end{multline}

For $\alpha<1/2$, the two-sided coverage for
$$
[\hat\theta-n^{-1/2}\hat\sigma z_{1-\alpha},\hat\theta + n^{-1/2}\hat\sigma z_{1-\alpha}]
$$
with $z_{1-\alpha}=-z_\alpha$ is $C(1-\alpha)-C(\alpha)$ which equals
\begin{multline*}
1-2\alpha 
+2n^{-1}z_{1-\alpha}\Bigl(\frac1{12}\kappa(z_{1-\alpha}^2-3) \\
-\frac1{18}\gamma^2\bigl( z_{1-\alpha}^4+2z_{1-\alpha}^2-3\bigr) 
-\frac14\bigl( z_{1-\alpha}^2+3\bigr)\Bigr)\varphi(z_{1-\alpha})+O(n^{-3/2}).  
\end{multline*} 
For the bootstrap methods covered by \eqref{eq:coverror}
the $O(n^{-3/2})$ error is actually $O(n^{-2})$ because the $O(n^{-3/2})$
terms cancel the same way that the $O(n^{-1/2})$ terms do.  This may very well be true
for `Norm' as well, but here we only consider terms up to order $O(1/n)$.

On page 949, \cite{hall:1988} writes:
``The term $-3.35$ appearing in the third column of
Table 1 arises from the difference between the standard normal distribution
function and an expansion of Student's t distribution function.''
This suggests that all we must do is drop the intercept from
the expression for the normal theory result, but for the avoidance
of doubt we work through the expansion here to verify that
this is the case.

The coverage error is the attained coverage minus the nominal
coverage.  For the two sided interval the coverage error is then
\begin{align}\label{eq:2sideerror}
\frac{2z_{1-\alpha}}n\Bigl(\frac{z_\alpha^2-3}{12}\kappa  
-\frac{ z_\alpha^4+2z_\alpha^2-3}{18} \gamma^2  
-\frac{z_\alpha^2+3}4
\Bigr)\varphi(z_\alpha)  
+O(n^{-3/2}) 
\end{align}  
where even powers of $z_{1-\alpha}$ and $z_\alpha$ are equal because $z_\alpha=-z_{1-\alpha}$.
The corresponding entry for Table~1, which leaves out the
factor $2\varphi(z_\alpha)/n$ is then
\begin{align}\label{eq:hallnorm}
\zoma\frac{z_\alpha^2-3}{12} \kappa
-\zoma\frac{z_\alpha^4+2z_\alpha^2-3}{18}\gamma^2-\zoma\frac{z_\alpha^2+3}2.
\end{align}
Using $z_{1-\alpha} = z_{0.975}\doteq 1.96$
we get $0.14\kappa -2.12\gamma^2-3.35$ after rounding,
matching Hall's Table 1.

For Gaussian data, we have $\gamma=\kappa=0$.  Then the standard normal
theory intervals~\eqref{eq:norm} have coverage error 
$$-\frac{2\varphi(\za)}n\zoma\frac{z_\alpha^2+3}{2}+O(n^{-3/2})
\doteq-\frac{0.78}n.$$
The negative sign indicates undercoverage.
The Student $t$ intervals~\eqref{eq:t} have coverage error zero.
Therefore we anticipate, and we will see, that the counterpart to
\eqref{eq:hallnorm} for Student $t$ intervals equals zero when $\gamma=\kappa=0$.

\section{Adjusted interval}\label{sec:adjusted}

To get the required formula,
we need to  replace $z_{1-\alpha}$ by $t_{(n-1)}^{1-\alpha}$,
the $1-\alpha$ quantile of Student's $t$ distribution.
We also have to replace $\hat\sigma$ by $s$.
We place $1-\alpha$ in the superscript of $t^{1-\alpha}_{(n-1)}$ so
it does not collide with the degrees of freedom $n-1$, while
placing $1-\alpha$ in ths subscript of $z_{1-\alpha}$ so that
it will not collide with exponents when we raise that quantile to integer powers.

Page 949 of \cite{abra:steg:1972} gives an asymptotic expansion
relating the $t$ quantiles to the Gaussian ones:
\begin{align}\label{eq:as72}
t_{(\nu)}^{\alpha} = z_\alpha 
+ \frac{g_1(z_\alpha)}\nu 
+ \frac{g_2(z_\alpha)}{\nu^2 }
+ \frac{g_3(z_\alpha)}{\nu^3}
+ \frac{g_4(z_\alpha)}{\nu^4}
+\cdots 
\end{align}
where 
\begin{align*}
g_1(x) &= \frac{x^3+x}4,\\
g_2(x) &= \frac{5x^5+16x^3+3x}{96},\\
g_3(x) &= \frac{3x^7+19x^5+17x^3-15x}{384},\quad\text{and}\\
g_4(x) &= \frac{79x^9+776x^7+1482x^5-1920x^3-945x}{92160}.
\end{align*}

Using the first two terms of expansion~\eqref{eq:as72},
$$
t^\alpha_{(n-1)} = z_\alpha + \frac{z_\alpha^3+z_\alpha}{4(n-1)}+ O(n^{-2})
= z_\alpha + \frac{z_\alpha^3+z_\alpha}{4n}+ O(n^{-2}).
$$
Similarly
$$
s = \hat\sigma\sqrt{\frac{n}{n-1}} = 1 + \frac1{2n} +O(n^{-2}),
$$
and so
\begin{align*}
st^\alpha_{(n-1)}
 &= \hat\sigma z_\alpha \Bigl(1+\frac1{2n}+O\Bigl(\frac1{n^2}\Bigr) \Bigr)
\Bigl(1+\frac{z_\alpha^2+1}{4n}+O\Bigl(\frac1{n^2}\Bigr)\Bigr)\\
 &= \hat\sigma z_\alpha \Bigl(1+\frac{z_\alpha^2+3}{4n}+O\Bigl(\frac1{n^2}\Bigr)\Bigr)\\
 &= \hat\sigma z_{\alpha'} 
\end{align*}
where
\begin{align*}
\alpha'&=\Phi(z_{\alpha'})=\Phi\Bigl(z_\alpha+\frac{z_\alpha^3+3z_\alpha}{4n}+O\Bigl(\frac1{n^2}\Bigr)\Bigr)\\
&=\alpha + \frac{z_\alpha^3+3z_\alpha}{4n}\varphi(z_\alpha) +O\Bigl(\frac1{n^2}\Bigr).
\end{align*}

Using~\eqref{eq:2sideerror}, the coverage error of the two sided Student $t$ interval is
\begin{align*}
&\phantom{=}\  2(\alpha-\alpha')+\frac{z_{1-\alpha'}}n\Bigl(\frac{z_{\alpha'}^2-3}{6}\kappa  
-\frac{ z_{\alpha'}^4+2z_{\alpha'}^2-3}9 \gamma^2  
-\frac{z_{\alpha'}^2+3}2  
\Bigr)\varphi(z_{\alpha'})  
+O(n^{-3/2}) \\
&=2(\alpha-\alpha')+\frac{z_{1-\alpha}}n\Bigl(\frac{z_{\alpha}^2-3}{6}\kappa  
-\frac{ z_{\alpha}^4+2z_{\alpha}^2-3}9 \gamma^2  
-\frac{z_{\alpha}^2+3}2  
\Bigr)\varphi(z_{\alpha})  
+O(n^{-3/2}). 
\end{align*} 
Now
\begin{align*}
2(\alpha-\alpha') &
= -\frac{z_\alpha^3+3z_\alpha}{2n}+O(n^{-2})
= \frac{z_{1-\alpha}}n\Bigl(\frac{z_\alpha^2+3}{2n}\Bigr)+O(n^{-2}).
\end{align*} 
It follows that the coverage error in Student's $t$ interval is
$$\frac{z_{1-\alpha}}n\Bigl(\frac{z_{\alpha}^2-3}{6}\kappa  
-\frac{ z_{\alpha}^4+2z_{\alpha}^2-3}9 \gamma^2\Bigr)\varphi(z_{\alpha})  
+O(n^{-3/2}). 
$$
The coefficients of $\kappa$ and $\gamma^2$  matches 
those of equation~\eqref{eq:2sideerror} for the normal theory
method and the intercept term is zero.

\section*{Acknowledgments}

This work was supported by the U.S.\ National Science Foundation
under grant DMS-2152780.

\bibliographystyle{apalike}
\bibliography{qmc}
\end{document}